\documentclass{amsart}

\usepackage[dvips]{graphicx}
\usepackage{amscd}
\usepackage{amsmath}
\usepackage{amsfonts}
\usepackage{amssymb}

\usepackage[all]{xy}

\def\openone
{\mathchoice
{\hbox{\upshape \small1\kern-3.3pt\normalsize1}}
{\hbox{\upshape \small1\kern-3.3pt\normalsize1}}
{\hbox{\upshape \tiny1\kern-2.3pt\SMALL1}}
{\hbox{\upshape \Tiny1\kern-2pt\tiny1}}}

\newtheorem{theorem}{Theorem}[section]
\newtheorem{corollary}[theorem]{Corollary}
\newtheorem{lemma}[theorem]{Lemma}
\newtheorem{proposition}[theorem]{Proposition}
\theoremstyle{definition}
\newtheorem{definition}[theorem]{Definition}
\newtheorem{remark}[theorem]{Remark}

\newtheorem{example}[theorem]{Example}

\theoremstyle{remark}

\def\1{{\openone}}
\def\A{{\mathfrak{A}}}
\def\R{{\bf R}}
\def\N{{\bf N}}
\def\Q{{\bf Q}}
\def\T{{\bf T}}
\def\Z{{\bf Z}}
\def\C{{\bf C}}
\def\ker{\mathop{\rm ker}\nolimits}
\def\tr{\mathop{\rm tr}\nolimits}
\def\id{\mathop{\rm id}\nolimits}
\def\Ad{\mathop{\rm Ad}\nolimits}
\def\Im{\mathop{\rm Im}\nolimits}

\title[Trace scaling automorphisms of certain stable AF algebras II]{Trace
scaling automorphisms of certain stable AF algebras II}
\author{Ola~Bratteli}
\address{Mathematics Institute\\
University of Oslo\\
PB 1053 Blindern\\
N-0316 Oslo\\ Norway}
\author{Akitaka~Kishimoto}
\address{Department of Mathematics\\
Hokkaido University\\
Sapporo\\ 060 Japan}

\begin{document}

\maketitle

\begin{abstract}
Two automorphisms of a simple stable AF-algebra with a
        finite dimensional lattice of lower semicontinuous traces are shown
to be
        outer conjugate if they act in the same way on $K_{0}$ and the
extremal traces are scaled by numbers which are not equal to 1 and satisfy a
certain condition (which holds if the scaling
        factors are all less than 1). The proof
        goes via the Rohlin property. As an application we consider the
        problem of classifying conjugacy or cocycle conjugacy classes of
        certain actions of $\T$ on a separable simple purely infinite
C*-algebra.
\end{abstract}

\section{Introduction}

This paper is a continuation of \cite{EEK} and \cite{EK}. In the
latter of these papers the case that $\A$ is a stable AF algebra with
totally ordered dimension group $K_{0}(\A)$ is considered. In this
case $\A$ has a densely defined lower semicontinuous trace $\tau$
(unique up to constant multiples) such that $K_{0}(\A)$ identifies
with the range $\tau_{\ast}(K_{0}(\A)) \subset \R$ , and in particular
$K_{0}(\A)$ has no infinitesimal elements (i.e., $\ker \tau_{\ast} =
0$). If $\alpha$ is an automorphism of $\A$, there exists necessarily
a $\lambda > 0$ such that $\tau_{\ast}\alpha_{\ast} = \lambda\tau_{\ast}$,
and it is proved in \cite{EK} that if $\lambda \not= 1$ then $\alpha$
has the following Rohlin property: For any $k \in \N$, any projection
$e \in \A$, any finite subset ${\mathcal F}$ of $\A$, and $\varepsilon > 0$,
there exists an orthogonal family $\{e_{ij}; i=1,2; j=0,1,\ldots=
 k_i-1\}$ of
projections in $\A$ with $k_{1} = k$, $k_{2} = k + 1$ such that
\begin{eqnarray*}
&&\sum_{i=1}^2\sum_{j=0}^{k_{i}-1} e_{ij} \geq e\;, \\ [.5ex]
&&\|\alpha (e_{ij}e)-e_{ij+1}\alpha(e)\|<
      \varepsilon\;,\qquad j= 0,1,\ldots , k_{i} - 2\;, \\ [.5ex]
&&\|\alpha ((e_{1 k_{1}- 1} + e_{2 k_{2}- 1}) e) - (e_{10} +
      e_{20}) \alpha (e) \| < \varepsilon\;, \\ [.5ex]
&&\|[x, e_{ij}] \| < \varepsilon\;,\qquad  x\in {\mathcal F}\;,
\end{eqnarray*}
where $e_{ik_{i}} = e_{i0}$. Here we may further impose the condition
$[e_{ij}, e] = 0$.

(This is slightly different from the Rohlin property given in 2.1 of
\cite{EK},  but this version easily follows and is strong enough to prove
3.1 of \cite{EK},  i.e., the stability of $\alpha$. See also \cite{K98},
where a relaxation  of this kind is discussed.)

It was also established in \cite{EK} that if $\A$ is any AF algebra
and $\alpha, \beta$ are automorphims of $\A$ with the Rohlin property,
then if $\alpha_{\ast} = \beta_{\ast}$ on $K_{0}(\A)$, $\alpha$ and
$\beta$ are outer conjugate: For any $\varepsilon > 0$ there is an
automorphism $\sigma$ of $\A$ and a unitary $U$ in $\A + \1$ such
that $\alpha = \Ad U \circ \sigma \circ \beta \circ \sigma^{-1}$,
$\|U-\1 \| < \varepsilon$ and $\sigma_{\ast}= \id$, so the previous
theorem  has an immediate corollary. The aim of the present paper is to
prove
the Rohlin property for trace scaling automorphisms under less
stringent assumptions on $\A$. Instead of assuming that $\A$ has one
trace separating elements in $K_{0}(\A)$ we will assume that the set
of lower semicontinuous traces on $\A$ form a finite dimensional lattice,
i.e. if $E$ is a fixed projection in the simple stable AF-algebra $\A$
then $\{ \tau | \tau$ is a trace on $\A$ and $\tau(E) = 1\}$ is a
finite-dimensional simplex. In order to prove the Rohlin property of
an automorphism $\alpha$, we must assume a scaling property of the
action $\alpha_{\ast}$ of $\alpha$ on the lattice of traces: As
$\alpha_{\ast}$ maps extremal traces into extremal traces, $\alpha_{\ast}$
permutes the finite number of extremal rays in the space of traces,
and hence some power of $\alpha_{\ast}$ scales all the extremal traces.
We choose the smallest positive integer $p$ such that $\alpha_{\ast}^p$
has this property and let $\Lambda$ be the set of scales;
$$
\Lambda = \{\lambda \mid \alpha_{\ast}^p (\tau) = \tau\alpha^p = \lambda\tau
\;\,\mbox{for an extreme trace $\tau \}$}\;.
$$
This is a finite subset of $(0,\infty)$.

\begin{definition}
        If $\A$ is a stable simple AF-algebra such that the set of lower
        semicontinuous traces form a finite dimensional lattice, and $\alpha$
        is an automorphism of $\A$, we say that $\alpha$ {\em scales the
        traces\/} if no extremal trace $\tau$ is invariant under any nonzero
        power $\alpha^n$ of $\alpha$. The set $\Lambda$ of positive numbers
introduced
        prior to the definition is called the {\em set of scales\/}
        for $\alpha$, and thus $\alpha$ scales the traces if and only if the
        set of scales does not contain 1.
\end{definition}

        To prove the main result in this paper we need more than that
        $\alpha$ scales the traces. The result is:

        \begin{theorem}
        Let $\A$ be a simple stable AF algebra such that the (densely
        defined) lower semicontinuos traces form a finite dimensional
        lattice. Assume that $\alpha$ scales the traces, and moreover that
        the ring $\Z[x,x^{-1}]$ of polynominals in $x$ and $x^{-1}$ with
        integer coefficients is dense in the algebra $C(\Lambda)$ of real
        functions on the scales. It follows that $\alpha$ has the Rohlin
        property.
        \end{theorem}

This theorem will be proved in Section~2.

        \begin{remark}
        The last condition in the theorem automatically implies that
        $1\notin\Lambda$, i.e. that $\alpha$ scales the traces. The
        condition is fulfilled if there is an element $p \in \Z[x,x^{-1}]$
        such that $0<|p(\lambda)| < 1$ for all $\lambda \in \Lambda $ by
        corollary 9.3 in \cite{Fer}, and such a $p$ automatically exist if in
        addition to $1 \notin \Lambda$ one has $\Lambda \subset (0,2)$,
        $\Lambda\subset ({1\over 2},\infty)$, or $\Lambda \subset \Q$.
        There is a finite set $\Lambda$ for which this condition is not
        fulfilled. See 2.8--9 below and \cite{Fer} for details.
        \end{remark}

We show in Example~2.10 below how to construct C*-dynamical systems
satisfying the hypotheses of Theorem~1.2 for a given $\Lambda$.

        \begin{corollary}
        Let $\A$ be a simple stable AF algebra such that the lower
        semicontinous traces form a finite dimensional lattice. Assume that
        the set $\Lambda$ of scales of $\alpha$ satisfies the condition in
        theorem~1.2. Let $\beta$ be another automorphism of $\A$ such that
        $\beta_{\ast} = \alpha_{\ast}$ on $K_{0}(\A)$. Then for any
        $\varepsilon > 0$ there is an automorphism $\sigma$ of $\A$ and a
        unitary $U$ in $\A + \1$ such that $\alpha = \Ad U \circ
        \sigma\circ \beta\circ\sigma^{-1}$, $\| U -\1\| < \varepsilon$ and
        $\sigma_{\ast} = \id$.
        \end{corollary}

\begin{proof}
This follows from Theorem~1.2 and \cite[Theorem~4.1]{EK}.
\end{proof}

        By combining this with a result of Krieger \cite{Kri}, we obtain a
        criterion for outer conjugacy of automorpisms of stationary systems
        \cite{Eff}, \cite{Tor}. To define these, let $\varphi$ be a $r
\times
        r$ matrix with nonnegative integer matrix elements, and let $G
(\varphi)$
        be the ordered abelian group defined as the inductive limit of
        $$
        \Z^r \stackrel{\varphi}{\longrightarrow} \Z^r
\stackrel{\varphi}{\longrightarrow} \Z^r \to
        $$
        where each $\Z^r$ is ordered by requiring non-negativity of the $r$
        coordinates. There is a unique stable AF algebra $\A_\varphi$
        associated to the dimension group $G(\varphi)$, and if $\varphi$ is
        primitive in the sense that $\varphi^n$ has strictly positive
        matrix elements for some $n \in \N$, the algebra $\A_{\varphi}$ is
        simple with a trace $\tau_{\varphi}$ which is unique up to a scale.
        In fact, if $\lambda_{\varphi}$ is the Perron-Frobenius eigenvalue
        of $\varphi$ and $\eta_{\varphi}$ is a corresponding left
        eigenvector,
$$
\eta_{\varphi}\varphi = \lambda_{\varphi}\eta_{\varphi}
$$
the trace $\tau_{\varphi}$ is given on the $n$'th group $\Z^r$ as
$$
        g \to \lambda_{\varphi}^{-n}\langle \eta_{\varphi}| g \rangle
$$
        Now, let $\sigma_{\varphi^\ast}$ be the natural shift automorphism of
        $G(\varphi)$ determined by the inductive limit diagram (see
        \cite[p.~37]{Eff}). Then
        $$
        \tau_{\varphi^\ast}\circ\sigma_{\varphi^\ast} =
        \lambda_{\varphi}^{-1}\tau_{\varphi^\ast}
        $$
        Unless $r = 1$ and $\varphi = 1$, we have $\lambda_{\varphi} > 1$.
        Let $\sigma_{\varphi}$ be some automorphism of $\A_{\varphi}$ such
        that $\sigma_{\varphi^\ast}$ is the shift automorphism of
$G(\varphi)$ defined above. For completeness we incorporate
        Krieger's result in the following corollary.

        \begin{corollary}
        Let $\varphi, \psi$ be primitive square matrics of dimension $r,s$,
respectively with
        non-negative integer matrix units and assume $(\varphi, r) \not=
        (1,1)$. Let $G(\varphi), G(\psi)$ be the assosiated dimension groups
        $G(\varphi),G(\psi)$, respectively, and let $\sigma_{\varphi},
\sigma_{\psi}$ be associated
        automorphisms of the corresponding stable AF algebras
        $\A_{\varphi}, \A_{\psi}$. Then the following conditions are
        equivalent:
\begin{itemize}
\item[(1)]
$\varphi$ and $\psi$ are shift equivalent, i.e. there exists $k \in \N$ and
there exists non-negative rectangular $r \times s$
        (resp. $s \times r$) matrices $A,B$ such that
\begin{eqnarray*}
&&A\psi = \varphi A , B\varphi = \psi B \\
&&AB = \varphi^k , BA = \psi^k
\end{eqnarray*}
\item[(2)]
$(G(\varphi), \sigma_{\varphi^\ast})$ and $(G(\psi),
        \sigma_{\psi^\ast})$
        are isomorphic, i.e. there exists an order isomorphism $\gamma_{\ast}:
        G(\psi) \to G (\varphi)$ with $\gamma_{\ast}\circ\sigma_{\psi*} =
        \sigma_{\varphi^\ast}\circ \gamma_{\ast}$

\item[(3)]
The C*-dynamical systems $(\A_{\psi}, \sigma_{\psi})$ and
        $(\A_{\varphi}, \sigma_{\varphi})$ are outer conjugate in the
        sense that there exists an isomorphism $\gamma: \A_{\psi}\to
        \A_{\varphi}$ and a unitary $U$ in $\A_{\varphi} + \1$ such that
        $\sigma_{\varphi} = \Ad U \circ \gamma \circ  \sigma_{\psi}\circ
        \gamma^{-1}$.
\end{itemize}
\end{corollary}

\begin{proof}
(1) $\Leftrightarrow$ (2) is in \cite[\S4.2]{Kri}, see also
\cite[Theorem~6.4]{Eff}.

(2) $\Rightarrow$ (3) follows from Corollary 1.4,
and (3) $\Rightarrow$ (2) is
        trivial, by defining $\gamma_{\ast}$ from $\gamma$.
        \end{proof}

        \begin{remark}
        The assumption $(\varphi,r) \not= (1,1)$ is irrelevant for the
        equivalence $(1) \Leftrightarrow (2)$. It is inserted to assume
that the
        Perron-Frahenius eigenvalue is larger than 1, such that
        $\sigma_{\varphi}$ really scales the trace by a scale less than 1.
        Note that $(2) \Rightarrow (3)$ actually is false in this
``trivial'' case,
        where $G = \Z$ and $\sigma_{\varphi^\ast} = \iota$ thus
$\A_{\varphi}   =
        \A_{\psi} =$ the algebra of compact operators. If we take
        $\sigma_{\psi} = \id$ and $\sigma_{\varphi} = \Ad(V)$, where $V-\1$ is
non-compact, it is clear that $\sigma_{\varphi}$ and
        $\sigma_{\psi}$ is not related as in (3) (Of course the Rohlin
        property also fails).
        \end{remark}

        \begin{remark}
        Note more generally that the condition of finite-dimensionality of
        the lattice of traces always is fulfilled if $\A =
        \overline{\bigcup_{n}\A_{n}}$ is an AF algebra, and the
        dimension of the center of $\A_{n}$ in the generating sequence is
        bounded in $n$ by a constant $K$. The projective dimension of the
        lattice of traces is then at most $K-1$. This is because the
        restriction of a trace to any $\A_{n}$ is a trace.

        In the situation of Theorem~1.2 it follows from \cite{R2} that the
        crossed product $\A \times \Z$ is purely infinite. As an application
        of Theorem~1.2 we will in Section~4 discuss conjugacy and cocycle
        conjugacy classes of certain actions of $\T$ on a separable simple
        purely infinite C*-algebra (whose crossed product is AF).

        Let us end the introduction by mentioning that a result analoguos to
        Corollary 1.4 has been proved by H.~Nakamura recently for the
        infinite case \cite{Nak}; If two automorpisms of a nuclear
separable
        simple purely infinite C*-algebra defines the same $KK$ class
        and have the Rohlin property, they are outer conjugate. (A similar
        result in the $A\T$ case is also found in \cite{K98}.)
        The expected characterization for the Rohlin property has been shown
        by M. Izumi \cite{Nak} i.e., if $\A$ is a separable purely infinite
        simple C*-algebra, and $\alpha$ is an automorphism of $\A$ such
        that $\alpha^n$ is outer for $n \not= 0$, then $\alpha$ has the
        Rohlin property. The proof goes roughly by lifting $\alpha$ to the
        C*-algebra of central sequences modulo trivial central
        sequences, which is known to be purely infinite \cite{KP}, and then
        prove the Rohlin property there without worrying about centrality.
        \end{remark}

\section{Proof of Theorem 1.2}

If $\alpha_{\ast}$ denotes the action of $\alpha$ on $K_{0}(\A)$, then
$K_{0}(\A)$ is a $\Z [x, x^{-1}]$ module as follows: If
$p(x,x^{-1}) \in \Z[x, x^{-1}]$ and $g \in K_{0}(\A)$, define the
action of $p$ on $g$ by
$$
pg = p(\alpha_{\ast}, \alpha_{\ast}^{-1})g\;.
$$
In the proof of Theorem~1.2 we will use the fact that $\Z[x,
x^{-1}]g$  is sufficiently big for each $g \in K_{0}(\A)_+\setminus\{0\}$.
The  condition on the scaling factors $\Lambda$ is imposed for this
purpose (cf. 2.7).

In proving Theorem~1.2 we will closely follow the method in Section 2
of \cite{EK} and recall the Rohlin property as defined there:

For any $k \in \N$ there are positive integers $k_1,\ldots, k_{m}\ge
k$ such that for any projections $E, e$ in $\A$, any unitary $U \in
\A + \C\1$, any finite subset ${\mathcal F}$ in $\A_{E}=E \A E$ and any
$\varepsilon > 0$ with $e \le E$, $\Ad U \circ \alpha (e) \le E$,
$e \in {\mathcal F}$, $\Ad U \circ \alpha (e) \in {\mathcal F}$ there exists
a family $\{e_{i,j} \mid i=1, \ldots, m,j=0, \ldots, k_i-1\}$ of
projections  in $\A$ such that
\begin{eqnarray*}
&&\qquad \sum\limits_{i}\sum\limits_{j}e_{i,j} = E \\
&&\| \Ad U \circ \alpha (e_{i,j}e) - e_{i, j+1} \Ad U \circ \alpha
(e) \| < \varepsilon \\
&&\qquad \|[x, e_{i,j}]|| < \varepsilon
\end{eqnarray*}
for $i = 1, \ldots, m$, $j = 0, \ldots , k_{i} -1$ and
$x \in {\mathcal F}$,
where $e_{i,k_i} = e_{i,0}$.
(Here we
can take $m=2$, $k_{1}=k$, and $k_{2}=k+1$.) Note that we have
imposed the condition
$$
\| \Ad U \circ \alpha (e_{i,k_i-1} e) - e_{i,0} \Ad U \circ
\alpha (e) \| < \varepsilon
$$
for each $i$, which is stronger than the corresponding part of the
Rohlin property introduced in Section 1:
$$
\| \Ad U \circ \alpha ((\sum_{i} e_{i,k_i-1})e) -
(\sum\limits_{i}e_{i,0}) \Ad U \circ \alpha (e) \| < \varepsilon\;.
$$
Our assumption on the set of scales will make it possible to prove
the Rohlin property in this stronger form. (We anticipate that the
weaker form of the Rohlin property stated in Section~1 could be true
with the special assumption on the set $\Lambda$ of scales in Theorem~2.1.
replaced by just $1\notin \Lambda$. This weaker form of the  Rohlin property
suffices for the proof of Corollary~1.4.)

By making an arbitrarily small inner perturbation of $\alpha$ and
choosing an increasing sequence $\A_{n}$ of finite-dimensional
subalgebras as a subsequence of any given such sequence, we may
assume that
\begin{eqnarray*}
&&\alpha^{-1}(\A_{n}) \subseteq \A_{n+1}\\
&&\alpha (\A_{n}) \subseteq \A_{n+1}
\end{eqnarray*}
for all $n$ (and $\bigcup_{n}\A_{n}$ is dense in $\A$).
Fix some non-zero projection $E$ in $\A_1$.
For given $n = 1,2, \ldots$, let $\A_{n} = \bigoplus_{j}\A_{n,j}$
be the central decomposition of $\A_{n}$ where $\A_{n,j}$ is a full
matrix algebra. Let $p_{j}^{(n)}\in \A_{n}\cap \A_{n}'$ be the
central projection corresponding to $\A_{n,j}$. By simplicity of $\A$,
we have either $p_{j}^{(n)} E = 0$ (this only happens when $n$ is
small), or
$$
K_{0}(\A p_{j}^{(n)} E) \simeq K_{0}(\A)
$$
Since $\alpha(\A_{n}) \subseteq \A_{n+1} \ , \
\alpha^{-1}(\A_{n}) \subseteq \A_{n+1}$
and
$$
[\alpha(x), y] = \alpha([x,\alpha^{-1}y])
$$
we have
$$
\alpha(\A \cap \A'_{n}) \subseteq \A \cap \A'_{n-1}
$$
for all $n$.

Now pick a projection $e$ in some $\A_{n}$ such that $e \le E$ and
$[\alpha(e)]\leq [E]$ and the central support of $e$ in $\A_{n E}$ is $E$.
By replacing
$n$ by a larger $n$  and modifying $\alpha$ by an
inner automorphism $\Ad U$ with $U \in\bigcup_{k}\A_{k} + \C\1$ we may
also assume that $\Ad U \circ \alpha (e) \le E$ and the central
support of $\Ad U \circ \alpha (e)$ in $\A_{n E}$ is $E$. Replacing
$\alpha$ by $\Ad U \circ \alpha$ and relabeling the algebras $\A_{n}$
we may thus assume
\begin{eqnarray*}
&&\quad e,\alpha(e), E \in \A_1\;, \\
&&e \le E\;\quad \alpha (e) \le E\;,
\end{eqnarray*}

\centerline{Central support of $e$ and $\alpha(e)$ in $\A_{1 E}$ (and thus
in all
$\A_{n E})$ is $E$.}
\bigskip

We now define a unital endomorphism.
$$
\widetilde{\alpha} : E (\A \cap \A_{n}')E \to E(\A \cap \A_{n-1}')E
$$
as follows: Because of simplicity of $\A$ and the central support
properties of $e , \alpha(e)$ we have
$$
E(\A \cap \A_{n}')E \cong e (\A \cap \A_{n}')e
$$
and
$$
E(\A \cap \A_{n-1}')E \cong \alpha (e) (\A \cap \A_{n-1}')\alpha(e)\;,
\qquad n\geq 2
$$

\medskip
\noindent
and since $\alpha$ maps $e(\A\cap\A_n')e$ into
$\alpha(e)(\A\cap\A_{n-1}')\alpha(e)$, we get a morphism
$$
\widetilde{\alpha} : E(\A \cap \A_{n}') E \to E(\A \cap \A_{n-1}')E
$$
by the requirement that the following diagram is commutative:
$$
\begin{array}{ccc}
E(\A \cap \A_n')E &
\stackrel{\widetilde{\alpha}}{\longrightarrow} &
     E(\A \cap \A_{n-1}')E \\
\downarrow\cong & & \downarrow\cong \\
e(\A \cap \A_{n}')e & \stackrel{\alpha}{\longrightarrow} &
\alpha(e)(\A\cap\A_{n-1}')\alpha(e)\end{array}
$$
Let us explain the vertical isomorphisms, for example the one to the
left: As $e \in \A_{1 E}\subseteq \A_{n E}$ has central support $E$
in $\A_{nE}$, it follows that $e$ has nonzero product with each of
the minimal central projections $p_{j}^{(n)}E$ in $\A_{nE}$ in
$\A_{nE}$. Let $(e_{k,l}^{(n,j)})_{1\le k,l\le [n,j]}$ be a complete
set of matrix units for $(\A_{n,j})_{E}$. There is a projection of
$\A_{E}$ onto $(\A \cap \A_{n}')_{E}$ given by
$$
x{\longrightarrow} \sum_{j,k}e_{k,1}^{(n,j)} x  e_{1,k}^{(n,j)}
$$
and as $e_{1,k}^{(n,j)}p_{i}^{(n)} = \delta_{i,j}e_{1,k}^{(n,j)}$ it
follows that
\begin{eqnarray*}
(\A \cap \A_{n}')p_{j}^{(n)} E &=& (\A \cap \A_{n,j}')p_{j}^{(n)}E \\
&\cong& e_{1,1}^{(n,j)}\A e_{1,1}^{(n,j)}
\end{eqnarray*}
Since $e$ has central support $E$, we may choose the matrix units
$(e_{k,l}^{(n,j)})$ such that $e_{1,1}^{(n,j)} \le e p_{j}^{(n)}$ and
$e p_j^{(n)}$ is the sum of the $e_{i,i}^{(n,j)}$ it has nonzero
product with. Thus repeating the argument above with $E$ replaced by
$e$, we deduce
$$
(\A \cap \A_{n}')p_{j}^{(n)}e \cong
e_{1,1}^{(n,j)}\A e_{1,1}^{(n,j)}
$$
we thus obtain the left vertical isomorphism.
\bigskip

We shall prove below that $\widetilde{\alpha}$, as a partially
defined endomorphism of $\A_{E}$, has the usual Rohlin property.
\bigskip

Note that $\widetilde{\alpha}$ defines a morphism from $E(\A \cap
\A_{n}')p_{j}^{(n)}E$ into $E(\A \cap \A_{n-1}')E$ by restriction
and since $p_{i}^{(n-1)} \in(\A \cap \A_{n-1}')'$,
$\,\widetilde{\alpha}$ defines a morphism
$$
\widetilde{\alpha}(i,j) : E(\A \cap \A_{n}')p_{j}^{(n)}E {\longrightarrow}
E(\A \cap \A_{n-1}')p_{i}^{(n-1)}E
$$
by multiplying $\widetilde{\alpha} (x)$ with $p_{i}^{(n-1)}$.
Let
$$
\widetilde{\alpha}_{\ast}(i,j) : K_{0}((\A \cap
\A_{n}')p_{j}^{(n)}E) \to
K_{0}((\A\cap\A_{n-1}')p_{i}^{(n-1)}E)
$$
be the corresponding map of the $K_{0}$-groups.
Recall that
$$
(\A \cap \A_{n}')p_{j}^{(n)}E = (\A \cap
\A_{n,j}')p_j^{(n)}E \cong  e_{1,1}^{(n,j)} \A e_{1,1}^{(n,j)}\;.
$$
Since $\A$ is a simple AF algebra, we deduce that
$$
K_{0}((\A \cap \A_{n}')p_{j}^{(n)}E) \cong K_{0}(e_{1,1}^{(n,j)}
\A e_{1,1}^{(n,j)})\cong K_{0}(\A)\;.
$$
Let us denote this isomorphism by $I_{\ast}^{(n,j)}$,
$$
I_{\ast}^{(n,j)}:K_{0}((\A \cap \A_{n}')p_{j}^{(n)}E)
\stackrel{\cong}{\longrightarrow}
K_{0}(\A)
$$
This should be distinguished from the monomorphism
$$
\iota_{\ast}^{(n,j)}:K_{0}((\A \cap \A_{n}')p_j^{(n)}E) \to
K_{0}(\A)
$$
which comes from the embedding
$$
(\A \cap \A_{n}')p_{j}^{(n)}E \subseteq \A
$$
In fact if $(\A_{n,j})_{E} = M_{[n,j]}$ then
$$
\iota_{\ast}^{(n,j)} = [n,j] I_{\ast}^{(n,j)}
$$
and $K_{0}((\A \cap \A_{n}')p_{j}^{(n)}E)$ identifies on the
one side with a subgroup of $K_{0}(\A)$ (through $\iota_{\ast}^{(n,j)}$)
and on the other side is isomorphic to $K_{0}(\A)$ (through
$I_{\ast}^{(n,j)}$).

Thus
$\widetilde{\alpha}_{\ast}(i,j): K_{0}((\A \cap
\A_{n}')p_j^{(n)}E) \to K_{0}((\A \cap
\A_{n-1}') p_{i}^{(n-1)}E)$
defines maps $\widehat{\alpha}_\ast(i,j):K_0(\A)\to K_0(\A)$ by
$\widehat{\alpha}_{\ast}(i,j)=I_\ast^{(n-1,i)}
\widetilde{\alpha}_{\ast}(i,j)(I_{\ast}^{(n,j)})^{-1}$.

\begin{lemma}
The map $\alpha_{\ast} : K_{0}(\A) \to K_{0}(\A)$ and the maps
$\widehat{\alpha}_{\ast}(i,j) : K_{0}(\A) \to K_{0}(\A)$ do all
commute:
$\widehat{\alpha}_{\ast}(i,j) \alpha_{\ast} =
\alpha_{\ast}\widehat{\alpha}_{\ast}(i,j)$.
\end{lemma}

\begin{proof}
We want to prove
$$
I_{\ast}^{(n-1,i)}\widetilde{\alpha}_{\ast}(i,j)(I_{\ast}^{(n,j)})^{-1}
\alpha_{\ast} =
\alpha_{\ast}I_{\ast}^{(n-1,i)}
\widetilde{\alpha}_{\ast}(i,j)(I_{\ast}^{(n,j)})^{-1}\;.
$$
Multiplying to the left by $(I^{(n-1,i)})^{-1}$ and to the right by
$I^{(n,j)}$, this amounts to showing
$\widetilde{\alpha}_{\ast}(i,j)((I_{\ast}^{(n,j)})^{-1}
\alpha_{\ast}I_{\ast}^{(n,j)}) =
((I_{\ast}^{(n-1,i)})^{-1}\alpha_{\ast}I_{\ast}^{(n-1,i)})
\widetilde{\alpha}_{\ast}(i,j)$ as a map from $K_{0}((\A \cap
\A_{n}')p_{j}^{(n)}E)$ into $K_{0}((\A \cap
\A_{n-1}')p_{i}^{(n-1)}E)$.
To this end, pick a
$$
g \in K_{0}((\A \cap \A_{n}')p_{j}^{(n)}E)
$$
such that
$$
0 \le g, \alpha_{\ast}^{(n,j)}(g) \le [p_{j}^{(n)}E]
$$
where $\alpha_\ast^{(n,j)}=(I_\ast^{(n,j)})^{-1}\alpha_\ast
I_\ast^{(n,j)}=(\iota_\ast^{(n,j)})^{-1}
\alpha_\ast \iota_\ast^{(n,j)}$ and $[\,\cdot\,]$
here refers to
$K_{0}$-class in
$(\A
\cap
\A_{n}')p_{j}^{(n)}E$. It suffices to verify the commutation on
these $g's$ since they generate $K_{0}((\A \cap
\A_{n}')p_{j}^{(n)}E)$.
Pick projections $h, f \in (\A \cap \A_{n}')p_{j}^{(n)}E$ such
that
$$
g = [h] \qquad \alpha_{\ast}^{(n,j)}(g) = [f]
$$
We may and will assume that
$h,f \in \bigcup\limits _{m=1}^\infty \A_{m}$. Note that
$$
[\alpha(h)]_{K_0(\A)}=\alpha_\ast \iota_\ast^{(n,j)}(g)
$$
and
\begin{eqnarray*}
[f]_{K_0(\A)} &=& \iota_\ast^{(n,j)}([f]) \\
&=&\iota_\ast^{(n,j)}\alpha_\ast^{(n,j)}(g) \\
&=&\alpha_\ast\iota_\ast^{(n,j)}(g)\;.
\end{eqnarray*}
If $\tau_{\ast}$ is any real-valued additive character on $K_{0}(\A)$ (no
positivity assumed), there is a unique linear functional $\tau$ on
$\bigcup\limits_{m=1}^\infty \A_{m}$ corresponding to $\tau_{\ast}$,
which restricts to a scalar multiple of the standard trace on each
$\A_{m,k}$, such that
$$
\tau_{\ast}([p]) = \tau (p)
$$
for each projection $p$ in $\bigcup\limits_{m=1}^\infty \A_m$.
To show our commutation, we just need to show
\begin{eqnarray*}
\lefteqn{\tau_\ast\iota_\ast^{(n-1,i)}
        \widetilde{\alpha}_\ast(i,j)([f])} \\
&&=\tau_\ast\iota_\ast^{(n-1,i)}\alpha_\ast^{(n-1,i)}
      \widetilde{\alpha}_\ast(i,j)(g) \\
&&=\tau_\ast\alpha_\ast\iota_\ast^{(n-1,i)}
     \widetilde{\alpha}_\ast(i,j)(g) \\
&&= \tau \circ \alpha(\widetilde{\alpha}(h)p_i^{(n-1)})
\end{eqnarray*}

\noindent
But this follows from the computation (where we use $\A_{n-1,i}$ to
denote $\A_{n-1,i}E)$
\begin{eqnarray*}
\lefteqn{\tau_{\ast}\iota_\ast^{(n-1,i)}
      (\widetilde{\alpha}_{\ast}(i,j)[f])} \\
&& =\tau(\widetilde{\alpha}(f) p_{i}^{(n-1)}) \\
&& ={\tau(\alpha(f e)p_{i}^{(n-1)})\over
     \tr_{\A_{n-1,i}}(\alpha(e))} \\
&& =\tau \circ \alpha (f)
    {\tr_{\A_{n,j}}(e \alpha^{-1}(p_{i}^{(n-1)}))\over
      \tr_{\A_{n-1,i}}(\alpha(e))} \\
&& =\tau \circ \alpha (\alpha(h))
     {\tr_{\A_{n,j}}(e \alpha^{-1}(p_{i}^{(n-1)}))\over
      \tr_{\A_{n-1,i}}(\alpha(e))} \\
&& =(\tau \circ \alpha) \circ \alpha(h)
     {\tr_{\A_{n,j}}(e \alpha^{-1}(p_{i}^{(n-1)}))\over
      \tr_{\A_{n-1,i}}(\alpha (e))} \\
&& ={\tau \circ \alpha (\alpha(he)p_{i}^{(n-1)})\over
      \tr_{\A_{n-1,i}}(\alpha (e))} \\
&& =\tau \circ \alpha (\widetilde{\alpha}(h)p_{i}^{(n-1)})
\end{eqnarray*}
Here the second  equality follows from the definition of
$\widetilde{\alpha}$, and the third from
$\tau (\alpha(f e) p_{i}^{(n-1)})
=\tau \circ \alpha (f e \alpha^{-1}(p_{i}^{(n-1)}))$
and the fact that $f\in(\A \cap \A_{n}')p_{j}^{(n)} E$ commutes
with $e$ and $\alpha^{-1}(p_{i}^{(n-1)})$, and $e$ commutes with
$\alpha^{-1}(p_{i}^{(n-1)}) \subseteq \A_{n-2}'$ so
$e\alpha^{-1}(p_{i}^{(n-1)})$ is a projection in a tensor product
complement of $f$.
\end{proof}

As in the case of $\alpha$ (but simpler) the embedding
$$
E((\A \cap \A_{n}')p_{j}^{(n)})E \subseteq E(\A \cap \A_{n-1}')E
$$
defines morphisms
$$
\iota_{\ast}(i,j): K_0((\A \cap \A_{n}')p_{j}^{(n)}E) \to K_0((\A \cap
\A_{n-1}')p_{i}^{(n-1)}E)
$$
by multiplying with $p_{i}^{(n-1)}$, and again $\iota_{\ast}(i,j)$ can
also be viewed as maps $K_{0}(\A) \to K_{0}(\A)$.

\begin{lemma}
The maps $\iota_{\ast}(i,j)$ commute with $\alpha_{\ast}$:
$$
\iota_{\ast}(i,j)\alpha_{\ast} = \alpha_{\ast}\iota_{\ast}(i,j)
$$
\end{lemma}

\begin{proof}
If $0\le g$, $\alpha_{\ast}g \le [p_j^{(n)}E]$, choose again
projections $h \in (\A \cap \A_{n,j}')p_j^{(n)}E$ and $f$ such
that
\begin{eqnarray*}
&&[h] = g \\
&&[f] = \alpha_\ast g
\end{eqnarray*}
Here we have identified $K_0((\A\cap\A_{n,j}')p_j^{(n)}E)$ with a subgroup
of
$K_0(\A)$ through $\iota_\ast^{(n,j)}$. Thus, for any real valued additive
character $\tau$ on $K_{0}(\A)$ we have again
\begin{eqnarray*}
\lefteqn{\tau (f p_{i}^{(n-1)}) = \tau(f)
     \tr_{\A_{n,j}}(p_{i}^{(n-1)})} \\
&& =\tau_\ast(\alpha_{\ast}(g))\tr_{\A_{n,j}}(p_{i}^{(n-1)}) \\
&& =\tau\circ\alpha(h)\tr_{\A_{n,j}}(p_{i}^{(n-1)}) \\
&& =\tau\alpha(h p_{i}^{(n-1)})
\end{eqnarray*}
establishing the lemma.
\end{proof}

\begin{lemma}
The morphisms
$$
\widetilde{\alpha}_\ast(i,j),\iota_\ast(i,j):
K_0((\A\cap\A_n')p_j^{(n)}E)\to K_0((\A\cap\A_{n-1}')p_i^{(n-1)}E)
$$
satisfy:
$$
\sum_{j}\widetilde{\alpha}_{\ast}(i,j)([p_{j}^{(n)}E]) =
[p_{i}^{(n-1)}E]
$$
and
$$
\sum_{j}\iota_{\ast}(i,j)
([p_{j}^{(n)}E]) = [p_{i}^{(n-1)}E]
$$
\end{lemma}

\begin{proof}
By definition of $\widetilde{\alpha} (i,j)$ we have
$$
\widetilde{\alpha}(i,j)(p_{j}^{(n)}E) = \widetilde{\alpha}
(p_{j}^{(n)}E)p_{i}^{(n-1)}E\;.
$$
But
\begin{eqnarray*}
\sum_j\widetilde{\alpha}(p_{j}^{(n)}E) &=& \widetilde{\alpha}
        (\sum_jp_{j}^{(n)}E) \\
&=& \widetilde{\alpha} (E) = E\;.
\end{eqnarray*}
so
\begin{eqnarray*}
\sum_j\widetilde{\alpha} (i,j) (p_{j}^{(n)}E) &=& E p_{i}^{(n-1)}E \\
&=& p_{i}^{(n-1)}E
\end{eqnarray*}
and hence
$$
\sum_{j}\widetilde{\alpha}_{\ast} (i,j)([p_{j}^{(n)}E]) =
[p_{i}^{(n-1)}E]
$$
Replacing $\widetilde{\alpha}$ by $\iota = \id$ in the reasoning above, we
obtain the other identity.
\end{proof}

Before continuing the proof we describe the action $\alpha_{\ast}$
on the lattice of traces in more details.

If $\R_{+}\tau_{k}$, $k=1,\ldots,d$ are the extreme rays of the
lattice of densely defined lower semi-continuous traces, then as
$\alpha_{\ast}$ maps extremal traces into extremal traces,
$\alpha_{\ast}$ permutes the extreme rays and we can choose the labels
$k$ such that $\alpha_{\ast}$ maps $\R_{+} \tau_{k}$ onto
$\R_{+}\tau_{k+1}$ except for a set $1\leq
k_{1}<k_{2} < \cdots < k_{l} = d$ such that
$\alpha_{\ast}(\R_{+}\tau_{k_{i}}) = \R_{+}\tau_{k_{i-1}+1}$ for
$i=1, \ldots , l$, where $k_{0} = 0$. But then
$\alpha_{\ast}^{k_i-k_{i-1}} (\tau_{k}) =
\lambda_{i}^{k_i-k_{i-1}} (\tau_{k})$ for $k_{i-1} < k \leq
k_{i}$ for a unique positive number $\lambda_{i}$.
If $p$ denotes the least common multiple of $k_i - k_{i-1}$,
$i=1,2, \ldots, l$, then $\alpha_{\ast}^p
(\tau_{k})$ is proportional to $\tau_{k}$ for all $k$ and $p$ is the
smallest positive integer with this property. Hence the set $\Lambda$ of
scales for
$\alpha$ is $\{\lambda_{i}^p : i=1,2,\ldots, l\}$.

To continue the proof we now normalize $\tau_{k}$ by $\tau_{k}(E)
=1$, where $E \in \A_{1}$ is fixed as before.

We view $\tau_{1}, \ldots, \tau_{d}$ as mutually disjoint extremal
trace states on $\A_{E}$, and the weak closures $\pi_{\tau_{i}}
(\A_{E})''$ are all isomorphic to the hyperfinite II$_1$ factor.
\bigskip

Let $\omega$ be a free ultrafilter on $\N$ and let $\A_{E \omega}$
be the C*-algebra of central sequences, i.e. $\A_{E \omega}$
consists of elements $(x_{n})$ in $l^\infty (\N, \A_{E})$ such
that $\lim\limits_{n\to \omega} \|x_{n}a - a x_{n}\| = 0$ for all
$a \in \A_{E}$ modulo sequences $(x_{n})$ such that $\lim\limits_{n \to
\omega} \| x_{n}\| = 0$. Since $\A = \overline{\bigcup_{n}\A_{n}}$ where
the $\A_n$ are finite dimensional, it follows that for any $x = (x_n) \in
\A_{E \omega}$ we can find an increasing sequence $k_{n}$ in $\N$ and
$x_{n}' \in (\A \cap \A_{k_n}')_{E}$ such that $k_{n}\to \infty$
and
$$
\lim_{n \to \omega} \|x_{n}-x_{n}'\| = 0
$$
Define a homomorphism $\widetilde{\alpha}_{\omega}$ of $\A_{E
\omega}$ into itself by
$$
\widetilde{\alpha}_{\omega}(x) = (\widetilde{\alpha} (x_{n}'))
$$
One now checks that $\widetilde{\alpha}_{\omega}$ is well defined in
the sense that it is independent of the earlier choice of $e$ and
$U$. Again following \cite{EK} one defines
$\widetilde{\alpha_{\omega}^{-1}}$ and verify
$\widetilde{\alpha_{\omega}^{-1}} \circ
\widetilde{\alpha}_{\omega} = \widetilde{\alpha}_{\omega}
\circ \widetilde{\alpha_{\omega}^{-1}} = \id$, so
$\widetilde{\alpha}_{\omega}$ is an automorphism of $\A_{E \omega}$.

Now define tracial states $\tau_{i\omega}$ on $\A_{E\omega}$ by
$$
\tau_{i \omega} (x) = \lim\limits_{n \to \omega} \tau_{i} (x_{n})
$$
for $x = (x_{n}) \in \A_{E \omega}$. We next show that
the $\tau_{i \omega}$ are permuted cyclically by
$\widetilde{\alpha}_{\omega}$, i.e.
$$
\tau_{i \omega} \circ \widetilde{\alpha}_{\omega} = \tau_{i+1 \omega}
$$
except for
$$
\tau_{k_i\omega} \circ \widetilde{\alpha}_{\omega} =
\tau_{k_{i-1}+1\omega}
$$
We first show that if $\tau$ is any extremal trace state on $\A_{E}$,
i.e. an extremal trace on $\A$ with $\tau(E) =1$, then
$$
\tau_{\omega} \circ \widetilde{\alpha} = (\tau
(\alpha(E)))^{-1}(\tau\circ\alpha)_{\omega}
$$
For this, pick any central sequence $x = (x_{n}) \in
\A_{E\omega}$ where $x_{n} \in (\A \cap \A_{k_{n}}')_{E}$
for some non-decreasing sequences $k_{n}$ with $k_{n}\to \infty$.
When  $k_{n} > k$ with $U\in\A_{k-1}+\C \1$ we have
$$
\tau (\widetilde{\alpha} (x_{n})) =
\sum_i \tau(\widetilde{\alpha}(x_{n})p_{i}^{(k-1)}E)
$$
Now
$$
y \in E (\A \cap \A_{k}')E \to {\tau (\widetilde{\alpha} (y)
p_{i}^{(k-1)}E) \over \tau (p_{i}^{(k-1)}E)}
$$
is a trace state. We have
$$
\widetilde{\alpha} (y) \in E (\A \cap \A_{k-1}')E
$$
so
$$
\widetilde{\alpha} (y) p_{i}^{(k-1)}E \in E (\A \cap
\A_{k-1}')p_{i}^{(k-1)}E\;.
$$
Now, $\widetilde{\alpha} (y) p_{i}^{(k-1)}E$ and $\Ad U\alpha(e)
p_i^{(k-1)}$ lie in different tensor factors of $\A_{E}p_i^{(k-1)}$ so
that
\begin{eqnarray*}
\lefteqn{\tau (\widetilde{\alpha}(y) \Ad U \alpha (e) p_{i}^{(k-1)})} \\
&&= \tau(\widetilde{\alpha}(y)p_{i}^{(k-1)})\tr_{\A_{k-1i,E}}
         (\Ad U\alpha(e)) \\
&&= \tau(\widetilde{\alpha}(y) p_{i}^{(k-1)})
{\tau(\Ad U \alpha(e) p_{i}^{(k-1)}) \over \tau (p_{i}^{(k-1)}E)}
\end{eqnarray*}
Thus, for large $n$,
$$
\tau (\widetilde{\alpha}(x_{n})) = \sum_i
{\tau(p_{i}^{(k-1)}E) \tau (\widetilde{\alpha}(x_{n}) \Ad U
\alpha (e)p_i^{(k-1)}E) \over
\tau (\Ad U \alpha (e) p_{i}^{(k-1)})}
$$
But
$$
\widetilde{\alpha} (x_{n}) \Ad U \alpha (e) = \Ad U \alpha
(x_{n} e)
$$
so, using the trace property:
$$
\tau (\widetilde{\alpha}(x_{n})) = \sum_{i}
{\tau (p_{i}^{(k-1)}E) \over \tau (p_{i}^{(k-1)}\alpha (e))}
(\tau \circ \alpha)(\alpha^{-1}(p_{i}^{(k-1)}) ex_n)
$$
Now, ${\tau \circ \alpha \over \tau (\alpha (E))}$ is a factor
state on $\A_E$, and hence
\begin{eqnarray*}
\lefteqn{\lim_{n\to\omega} {\tau\alpha (\alpha^{-1}(p_{i}^{(k-1)})
    ex_n)\over \tau (\alpha(E))}} \\
&&={\tau\alpha(\alpha^{-1}(p_{i}^{(k-1)})e) \over
\tau(\alpha(E))}
\cdot {(\tau\circ\alpha)_\omega(x) \over \tau (\alpha (E))}
\end{eqnarray*}
and hence, by the previous relation,
$\lim\limits_{n \to \omega} \tau (\widetilde{\alpha} (x_{n})) =
\sum\limits_{i} {\tau (p_{i}^{(k-1)}E) \over \tau(\alpha(E))}
(\tau\alpha)_{\omega}(x)
= {(\tau \alpha)_{\omega}(x) \over \tau(\alpha(E))}\,$,
or:
$$
\tau_{\omega}\circ \widetilde{\alpha}_{\omega} =
(\tau(\alpha(E)))^{-1}(\tau\circ\alpha)_{\omega}\;.
$$

It follows immediately that the trace states $\tau_{i \omega}$ are
permuted cyclically (without scaling) by
$\widetilde{\alpha}_{\omega}$.
The states $\tau_{i \omega}$ are trace states on $\A_{E \omega}$
because $\varphi_{\omega}$ is a trace state on $\A_{E \omega}$ for
{\em any\/} state $\varphi$ on $\A_{E}$, see \cite{C74}. Furthermore, if
$\varphi_{1}, \varphi_{2}$ are disjoint states on $\A_{E}$, the
traces $\varphi_{1 \omega}$ and $\varphi_{2 \omega}$ are disjoint
on $\A_{E \omega}$. To see this, let $F$ be the central support
for $\varphi_{1}$ in $\A_{E}^{\ast\ast}$, and pick a sequence $x_{n} \in
\A_{E}$ such that $\|x_n\| \le 1$ and $s-\lim\limits_{n \to
\omega} x_{n}= F$.
Since $F$ is central, we have $U F U^{\ast} = F$ for all $U
\in {\mathcal U} (\A_{mE})$, and by approximating the projection
$\int\limits_{{\mathcal U}(\A_{mE})} U\cdot U^{\ast} dU$ onto $\A_{mE}'$ by
finite convex combinations, we may  replace $x_{n}$ by a sequence such that
there is a sequence
$k_{n} \in \N$ such that $(k_{n})$ is increasing, $k_{n} \to \infty$ and
$x_{n} \in E(\A \cap \A_{k_n}')E$ for all $n$. Then the new sequence
$x=(x_{n})$ is norm central and corresponds to an element
$\widetilde{F}$ in the unit ball of $\A_{E\omega}$. But then
$\varphi_{1 \omega}(\widetilde{F}) = 1,\varphi_{2
\omega}(\widetilde{F}) = 0$, so the trace states
$\varphi_{1\omega}$ and $\varphi_{2\omega}$ are disjoint.
We have defined the tracial states $\tau_{1 \omega},\ldots ,
\tau_{d\omega}$ of $\A_{E\omega}$ and shown that they are mutually
disjoint. Let $\pi_{i\omega}$ be the $GNS$ representation of $\A_{E
\omega}$ associated to $\tau_{i\omega}$ and consider the sum
$\pi_{\omega}=\bigoplus\limits_{i=1}^d \pi_{i \omega}$. If $\,\mathfrak{R}
= \pi_{\omega} (\A_{E\omega})'' = \bigoplus\limits_{i=1}^d
\pi_{i \omega}(\A_{E\omega})''$, $\widetilde{\alpha}_{\omega}$
extends to an automorphism of $\,\mathfrak{R}$, which is denoted again by
$\widetilde{\alpha}_{\omega}$.

\begin{lemma}
For any $n \in \Z \setminus \{0\}$, $\widetilde{\alpha}_{\omega}^n$ is a
properly outer automorphism of $\,\mathfrak{R}$.
\end{lemma}

\begin{proof}
Suppose that $\pi_{i\omega}\circ \widetilde{\alpha}_{\omega}^n =
\pi_{i\omega}$.
Then as in the proof of Lemma 2.4 of \cite{EK} we argue that
$\widetilde{\alpha}_{\omega}^n | \pi_{i\omega}(\A_{E\omega})''$ is
properly outer (where we use the fact that
$\pi_{i\omega}(\A_{E\omega})''
=\pi_{i\omega}(\A_{E\omega})=\mathfrak{R}_{i\omega}$ for
$\,\mathfrak{R}_{i} = \pi_i(\A_{E})'')$. Thus we conclude that
$\widetilde{\alpha}_{\omega}^n$ is properly outer on $\,\mathfrak{R}$.
\end{proof}

\begin{lemma}
Recall that
$k_{1}, k_{2}-k_{1}, k_{3}-k_{2}, \cdots, k_{l} - k_{l-1}$
are the periods for the cycles in $\{\tau_{1 \omega},
\ldots, \tau_{d \omega}\}$ under the action of
$\widetilde{\alpha}_{\omega}$. Let again $p$ be the least common
multiple of $k_{1}, k_{2}-k_{1}, \ldots, k_{l}-k_{l-1}$. Then for
any $n \in \N$ there are $n p$ projections $F_{0},F_{1}, \ldots ,
F_{n p-1}$ in $\,\mathfrak{R}$ such that
\begin{eqnarray*}
&&\sum_{i=0}^{n p-1}F_{i} = \1\;, \\
&&\widetilde{\alpha}_{\omega}(F_{i})=F_{i+1}\;, \qquad
      i= 0,1,\ldots,np-1 \\
&&\tau_{k\omega}(F_{i}) ={1 \over np}\;, \qquad k =1,2,\ldots,d
\end{eqnarray*}
where $F_{np} = F_{0}$.
\end{lemma}

\begin{proof}
Without the condition $\tau_{k \omega}(F_{i}) = {1 \over n p}$,
this follows from \cite{C75}.

To obtain this condition we consider each cycle $\{\tau_{k_{i+1}
\omega}, \ldots , \tau_{k_{i+1}}\}$ separately, and thus we may suppose that
$\{\tau_{1 \omega}, \ldots , \tau_{\alpha \omega}\}$ is one cycle,
i.e., $\tau_{i \omega} \circ \widetilde{\alpha} = \tau_{i+1\omega}$
except $\tau_{d\omega}\circ \widetilde{\alpha}_{\omega} =
\tau_{1 \omega}$.
In this case we find $nd$ projections
$F_{i}^{(j)}$, $j = 1, 2, \ldots , d$, $i=0,1, \ldots, n - 1$
in $\pi_{1 \omega}(\A_{E \omega})''$ such that
\begin{eqnarray*}
&&\sum_{j=1}^d \sum_{i=0}^{n-1} F_{i}^{(j)} = \1 \qquad
        \mbox{in $\,\pi_{1\omega} (\A_{E\omega})''$}\;, \\
&&\widetilde{\alpha}_{\omega}^d (F_{i}^{(j)}) = F_{i+1}^{(j)}\;,\qquad
       i = 0,1, \ldots , n-1
\end{eqnarray*}
where $F_{n}^{(j)} = F_{0}^{(j)}$. Then we set
$$
F_{i} = \sum\limits_{j=1}^d \widetilde{\alpha}_{\omega}^{i+j-1}
(F_{0}^{(j)})
$$
for $i= 0, 1, \ldots, nd -1$. It follows that
\begin{eqnarray*}
&&\sum_{i=0}^{nd-1} F_{i} = \1 \qquad \mbox{in $\,\mathfrak{R}$}\;, \\
&&\widetilde{\alpha}_\omega (F_{i}) = F_{i+1}\;, \\
&&\tau_{k\omega} (F_{i}) = {1 \over nd}\;.
\end{eqnarray*}
\end{proof}

\begin{lemma}
There exists an orthogonal family $\{f_0,f_1,\ldots,f_{np-1}\}$ of
projections in $\A_{E\omega}$ such that
\begin{eqnarray*}
&&\widetilde{\alpha}_\omega(f_i)=f_{i+1}\qquad
       i=0,1,\ldots,np-2\;. \\
&&\tau_{k\omega}(f_i)=\frac{1}{np}\;.
\end{eqnarray*}
\end{lemma}

\begin{proof}
For each $F_i$ in Lemma~2.5 we find a projection
$f_i\in\A_{E\omega}$ such that $F_i=\pi_\omega(f_i)$. However this
only shows that
$\pi_\omega(\widetilde{\alpha}_\omega(f_i))=
\pi_\omega(f_{i+1})$, and $\pi_\omega$ will never be
faithful. To overcome this difficulty we replace $f_0$ by a smaller
projection and we set $f_i=\widetilde{\alpha}_\omega^i(f_0)$ as in
\cite[K96]{K95}. Note that we have not assumed that
$\widetilde{\alpha}_\omega(f_{np-1})=f_0$.
\end{proof}

\begin{lemma}
Recall the overall assumptiuon in Theorem~2.1 that the image of
$\Z[x,x^{-1}]$ in $C(\Lambda)=C$
$(\{\lambda_1^p,\ldots,\lambda_l^p\})$ is dense. There exists a
sequence $(f_{0m})$ of projections in $\A_{E}$ such that
$f_{0m}\in(\A_E\cap\A_{k_{m}}')_E$ with $k_m\to\infty$ as $m\to\infty$,
$(f_{0m})$ in $\A_{E\omega}$ is less than $f_0$ in Lemma~2.6,
$$
\tau_{k\omega}((f_{0m}))=\frac{1}{np}\qquad
\mbox{for $\,k=1,2,\ldots,d$}\;,
$$
and
$$
[f_{0m}p_j^{(k_m)}]= p_m(\alpha_\ast^p,\alpha_\ast^{-p})[p_j^{(k_m)}]\;,
$$
in $K_0((\A\cap\A_{k_m}')_E)$ for some $p_m\in\Z[x,x^{-1}]$, where
$p_j^{(k_m)}$ denotes $p_j^{(k_m)}E$.
\end{lemma}

\begin{proof}
First we take a sequence $(f_{0m})$ of projections in $\A_E$ which
represents $f_0\in\A_{E\omega}$ in Lemma~2.6 and satisfies that
$f_{0m}\in(\A\cap\A_{k_m}')_E$ and for $p_j^{(k)}$ we have that
$$
\tau_s(f_{0m}p_j^{(k)})\to\frac{1}{np} \tau_s(p_j^{(k)})\qquad
\mbox{as $\;m\to\infty$}\;.
$$
For any $\varepsilon>0$ with $\varepsilon<\frac{1}{2np}$ it follows from
the density of $\Z[x,x^{-1}]$ in $C(\Lambda)$ that
there exists a function $q\in \Z[x,x^{-1}]$ such that
$$
\frac{1}{np}-2\varepsilon<q(\lambda_i^p,\lambda_i^{-p})<
\frac{1}{np}-\varepsilon\;.
$$
If $k_{i-1}<s\leq k_i$ then
$$
\tau_s(q(\alpha^p,\alpha^{-p})p_j^{(n)})=
q(\lambda_i^p,\lambda_i^{-p})\tau_s(p_j^{(n)})
$$
since $\tau_s\alpha^p=\lambda_i^p\tau_s$. Thus for a sufficiently
large $m$, it follows that
$$
q(\alpha_\ast^p,\alpha_\ast^{-p})[p_j^{(n)}]\leq
[f_{0m}p_j^{(n)}]
$$
in $K_0((\A\cap\A_k')_E$. Then we find a subprojection $f_{0m}'$ of
$f_{0m}$ in $(\A\cap\A_n')_E$ such that
$$
[f_{0m}' p_j^{(n)}]=
q(\alpha_\ast^p,\alpha_\ast^{-p})[p_j^{(n)}]\;.
$$
For an increasing $k$ and a decreasing $\varepsilon$ we use the
above argument to replace $f_{0m}$ by a smaller projection
satisfying the required properties.
\end{proof}

\begin{proof}[Proof of Theorem~1.2]
By using the projections $f_{0m}$ in Lemma~2.7 we obtain that
$$
f_{0m},\widetilde{\alpha}(f_{0m}),\widetilde{\alpha}^2(f_{0m}),\ldots,
\widetilde{\alpha}^{np-1}(f_{0m})
$$
are almost mutually orthogonal when $m\to\omega$ and we must show that they
are equivalent as projections in $(\A\cap\A_{k_m-np}')_E$.

To show that they are equivalent, we regard $\widetilde{\alpha},\iota$
as homomorphisms of $(\A\cap\A_{k_m}')_E$ into $(\A\cap\A_{k_m-1}')_E$ and
compute by Lemmas~2.1--2.3
\begin{eqnarray*}
\widetilde{\alpha}_\ast([f_{0m}]) &=&
    \bigg( \sum_j\widetilde{\alpha}_\ast(i,j)
    ([f_{0m}p_j^{(k_m)}])\bigg)_i \\
&=& \bigg( \sum_j\widetilde{\alpha}_\ast(i,j)
   (p_m(\alpha_\ast^p,\alpha_\ast^{-p})[p_j^{(k_m)}])\bigg)_i \\
&=& \bigg( p_m(\alpha_\ast^p,\alpha_\ast^{-p})
     \sum_j\widetilde{\alpha}_\ast(i,j)
       [p_j^{(k_m)}])\bigg)_i \\
&=& (p_m(\alpha_\ast^p,\alpha_\ast^{-p})[p_i^{(k_m-1)}])_i\;.
\end{eqnarray*}
With a similar computation for $\iota_\ast$ we obtain that
\begin{eqnarray*}
[\widetilde{\alpha}(f_{0m})] &=&
                         \widetilde{\alpha}_\ast([f_{0m}]) \\
&=& (p_m(\alpha_\ast^p,\alpha_\ast^{-p})[p_i^{(k_m-1)}])_i \\
&=& \iota_\ast([f_{0m}])\;,
\end{eqnarray*}
i.e., $[\widetilde{\alpha}(f_{0m})]=[f_{0m}]$ in
$K_0((\A\cap\A_{k_m-1}')_E)$. We can repeat this process.
Thus we find a C*-subalgebra $\mathfrak{B}$ of $(\A\cap\A_{k_m-np}')_E$
such that $\mathfrak{B}$ is a factor of type $I_{np}$ which almost contains
the
$np$ projections
$f_{0m},\widetilde{\alpha}(f_{0m}),\ldots,
\widetilde{\alpha}^{np-1}(f_{0m})$, the identity of $\mathfrak{B}$ is
close to $E$ when the distance is measured with the traces $\tau_s$, and
the large portion of $\mathfrak{B}$ is $\widetilde{\alpha}$-invariant.
>From this one can deduce the Rohlin property for $\widetilde{\alpha}$ as
in
\cite[K96]{K95}. By reinterpreting $\widetilde{\alpha}$ in terms of $\Ad
U\circ \alpha$ and $e$, one obtains the Rohlin property of $\alpha$.
\end{proof}

After having finished the proof of Theorem~1.2 we now consider the special
condition on the ring $\Z[x,x^{-1}]$ of Laurent polynomials, i.e. we
establish the statements in Remark~1.3.

\begin{proposition}
Let $\Lambda$ be a compact subset of $(0,1)\cup(1,\infty)$. Then
$\Z[x,x^{-1}]$ is dense in $C(\Lambda)$ if one of the following
conditions hold:
\begin{enumerate}
\item
$\Lambda\subset C(0,2)$.
\item
$\Lambda\subset(1/2,\infty)$.
\item
$\Lambda$ is a finite subset of the rational numbers $\Q$.
\end{enumerate}
\end{proposition}

\begin{proof}
By Corollary~9.3 of \cite{Fer} it suffices to show that there is a
$p\in\Z[x,x^{-1}]$ such that $\|p\|_\Lambda<1$ and
$p(\lambda,\lambda^{-1})\not=0$ for $x\in\Lambda$.

(The existence of $p$ is clearly necessary for the density. For the
sufficiency, let us extract the apposite argument from \cite{Fer}:
Let $f\in C(\Lambda)$ and put $g=f/p$. Then $g\in C(\Lambda)$. By
Weierstrass theorem, there exists a real polynomial
$q(x)=\sum\limits_{k=0}^r \beta_kx^k$ approximating $g$, and we may
assume that each $\beta_k$ has the form $\beta=n2^{-m}$. If $a=\inf p>0$
 and
$b=\sup p<1$ we may approximate each $\beta$ arbitrarily well on $[a,b]$ by
$$
n\frac{1}{(2-t^k)^m}=
\frac{n}{(1-(t^k-1))^m}=
n\bigg( \sum_{i=1}^\infty(t^k-1)^i\bigg)^m
$$
by choosing $k$ large and replacing $\infty$ by a finite number, and
thus we may approximate the first order polynomials $\beta_kt$ arbitrarily
well by polynomials $p_k(t)\in\Z[t]$ for $t\in[a,b]$. Thus
$$
q(x)p(x)-\sum_{k=0}^r x^k p_k(p(x))=
\sum_{k=0}^r x^k(\beta_kp(x)-p_k(p(x)))
$$
is uniformly small for $x\in\Lambda$, and since
$f(x)=g(x)p(x)\approx q(x)p(x)$ and $\sum\limits_{k=0}^r x^k
p^k(p(x))\in\Z[x,x^{-1}]$, the sufficiency follows.)

We may use $p=x(2-x)$ to establish (1) and $p=\frac{1}{x}(2-\frac{1}{x})$
to establish (2) in the Proposition.

In case (3), let $\Lambda\cap(0,1)=\{\frac{p_i}{q_i}: i=1,\ldots,k\}$
and $\Lambda\cap(1,\infty)=\{\frac{p_i}{q_i}: i=k+1,\ldots,\ell\}$,
where $p_i,q_i$ are integers, and let
$$
P(x,x^{-1})=\frac{1}{x^n}\prod_{i=1}^k
(q_ix-p_i)^2+
x^n \prod_{i=k+1}^\ell(q_ix-p_i)^2
$$
for a sufficiently large $n$.
\end{proof}

\begin{example}
If $\Lambda=\{2-\sqrt{3},2+\sqrt{3}\,\}$ the image of $\Z[x,x^{-1}]$
in $C(\Lambda)$ is not dense.

Note that $\lambda_1=2-\sqrt{3}$ and
$\lambda_2=\frac{1}{\lambda_1}=2+\sqrt{3}$ are the roots of the
equations
$$
x^2-4x+1=0
$$
and
$$
\frac{1}{x^2}-\frac{4}{x}+1=0\;.
$$
Thus it suffices to show that the image of
$$
\{ax+\frac{b}{x}+c|a,b,c\in\Z\}
$$
in $C(\{\lambda_1,\lambda_2\})$ is not dense, or
$$
\{(a\lambda_1+b\lambda_2+c,\; a\lambda_2+b\lambda_1+c)\mid a,b,c\in\Z\}
$$
is not dense in $\R^2$. Suppose that for some $a,b,c\in\Z$.
$$
0<a\lambda_1+b\lambda_2+c<1,\quad
0<a\lambda_2+b\lambda_1+c<1\;.
$$
If $a=b$, then $a\lambda_1+b\lambda_2+c=4a+c\in\Z$, which is a
contradiction. Hence $a\not= b$. Since
$$
(a-b)\lambda_1-(a-b)\lambda_2=2(a-b)\sqrt{3}
$$
has modulus at least 1, we reach a contradition. Thus the assertion
follows.
\end{example}

\begin{example}
In this example we show how to construct C*-dynamical systems $(\A,\alpha)$
satisfying the hypotheses in Theorem~1.2. Let
$\Lambda$ be a finite subset of
$(0,1)\cup(1,\infty)$ and let
$d\in\N$. We define an order on $\bigoplus\limits_1^d\Z[x,x^{-1}]$ by
$p\geq q$ if $p_i(\lambda)>q_i(\lambda)$ for $i=1,2,\ldots,d$ and
$\lambda\in\Lambda$ or $p_i\equiv q_i$ for $i=1,2,\ldots,d$ where
$p=(p_1,p_2,\ldots,p_d)$ etc. If the image of $\Z[x,x^{-1}]$ in
$C(\Lambda)$ is dense, then this order on
$\bigoplus\limits_1^d\Z[x,x^{-1}]$ gives a dimension group, which
will be denoted by $G$. We define an automorphism $\alpha_\ast$ of
$G$ by
$$
\alpha_\ast((p_1,p_2\ldots,p_d))=(xp_d,p_1,p_2,\ldots,p_{d-1})\;.
$$
$G$ is associated to a stable AF algebra $\A$ and $\alpha_\ast$ to an
automorphism $\alpha$ of $\A$. For each $i=1,2,\ldots,d$ and
$\lambda\in\Lambda$ there is a positive character $\tau_{i,\lambda}$
on $G$ such that $\tau_{i,\lambda}(p)=p_i(\lambda)$ for
$p=(p_1,p_2,\ldots,p_d)\in G$. Each $\tau_{i,\lambda}$ gives a lower
semicontinuous trace on $\A$, which we denote by the same symbol.
Since
$$
\tau_{i,\lambda}\alpha=\left\{
\begin{array}{ll}
\tau_{i-1,\lambda} & \mbox{for $\,i=2,3,\ldots,d$} \\
\lambda\tau_{d,\lambda} & \mbox{for $\,i=1\;,$}\end{array} \right.
$$
it follows that $\tau_{i\lambda}\alpha^d=\lambda\tau_{i\lambda}$ for
all $i,\lambda$ and thus the set of scales is $\Lambda$. Since
$(\A,\alpha)$ satisfies the assumptions in Theorem~1.2, $\alpha$ has
the Rohlin property. In passing we can see that
$$
G/(\id-\alpha_\ast)G\cong \Z\;.
$$
Since $\A\times_\alpha\Z$ is purely infinite \cite{R1},
$\A\times_\alpha\Z$ must be isomorphic to $O_\infty\otimes \mathfrak{K}$
\cite{KP}. (This conclusion also holds if $\Lambda$ is a non-finite closed
subset of
$(0,1)\cup(1,\infty)$ such that the dimension group condition is
satisfied, as it is when $\bigoplus\limits_1^d\Z[x,x^{-1}]$ is dense in
$C(\Lambda)$.)
\end{example}

\begin{remark}
When there is an automorphism $\beta$ of $K_0(\A)$ such that
$\beta\alpha_\ast=\alpha_\ast\beta$, we obtain that
$\widetilde{\alpha}_\ast(i,j)\beta=
\beta\widetilde{\alpha}_\ast(i,j)$ and
$\iota_\ast(i,j)\beta=\beta \iota_\ast(i,j)$ as in Lemmas~2.1 and 2.2 and
we can sometimes use this fact in Lemma~2.7 to prove that $\alpha$ has the
Rohlin property. As an example, take a $\lambda\in(0,1/2)$ and equip the
abelian group $\Z[x,x^{-1},(1-x)^{-1}]$ with the order defined by
$p\geq 0$ if $p(\lambda,\lambda^{-1},(1-\lambda)^{-1})>0$ and
$p(\frac{1}{2},2,2)>0$, or $p\equiv 0$. Let $\A$ be the stable AF
algebra corresponding to this ordered group and $\alpha$ an
automorphism of $\A$ corresponding to the multiplication of
$\frac{1-x}{x}$ such that in the tracial representation
corresponding to $1/2$ all non-zero poweres of $\alpha$ are not
weakly inner. (For the tracial representation corresponding to $\lambda$
this follows automatically.) The automorphism $\beta$ of $K_0(\A)$
corresponding to multiplication by $x$ commutes with $\alpha_\ast$.
To prove that $\alpha$ has the Rohlin property we proceed exactly
as before but use polynomials in $\beta$ and $\beta^{-1}$ in
Lemma~2.7. (Since $\beta$ scales the traces by $1/2$ and $\lambda$,
the required condition for 2.7 is satisfied;, the image of
$\Z[x,x^{-1}]$ in $C(\{\lambda,1/2\})$ is dense. In passing we note
that the quotient of $K_0(\A)$ by $\Im(\id-\alpha_\ast)$ is
isomorphic to $\Z[\frac{1}{2}]$, the dimension group of the UHF
algebra of type $2^\infty$. Since $\ker(\id-\alpha_\ast)=0$, it
follows that $\A\times_\alpha\Z$ has the same $K$-theory as the UHF
algebra of type $2^\infty$. (The automorphism $\alpha$ was first
used in \cite{BEH}.)
\end{remark}

\begin{remark}
By taking the tensor product of $\A$ with a UHF algebra $\mathfrak{B}$, we
could overcome the difficulty encountered in Lemma~2.7. We could
thus show in Theorem~1.2 without the condition on the scaling
factors that $\alpha\otimes\id$ on $\A\otimes \mathfrak{B}$ has the Rohlin
property.
\end{remark}

\section{Actions of a compact abelian group}

We will describe here some invariants for conjugacy and cocycle
conjugacy classes of actions of a compact abelian group on a unital
C*-algebra, borrowing ideas from \cite{Bla}, \cite{HR}.

Let $\A$ be a unital C*-algebra and $\alpha$ an action of
a compact abelian group $G$ on $\A$. The dual system is
$(\A\times_\alpha G,\widehat{G},\widehat{\alpha})$, where
$\widehat{G}$ is the dual group of $G$. The dual system of
$(\A\times_\alpha G, \widehat{G},\widehat{\alpha})$ is
isomorphic to $(\A\otimes \mathfrak{K}(L^2(G)), G,
\alpha\otimes\Ad \widetilde{\lambda})$ by Takai's duality \cite{Tak},
where
$\mathfrak{K}(L^2(G))$ is the compact operators on $L^2(G)$ and
$\widetilde{\lambda}$ is the unitary representation of $G$ on $L^2(G)$
defined by
$$
(\widetilde{\lambda}_t\xi)(s)=\xi(t+s)\;.
$$
>From the inclusion $\iota:\A\times_\alpha G\subset\A\otimes
\mathfrak{K}(L^2(G))$ we obtain the natural map $\iota_\ast:
K_i(\A\times_\alpha G)\to K_i(\A)$ for
$i=0,1$ which satisfies that $\iota_\ast
\widehat{\alpha}(s)_\ast=\iota_\ast$, $s\in\widehat{G}$.

Let $\lambda$ be the canonical unitary representation of $G$ in the
multiplier algebra $M(\A\times_\alpha G)$ of
$\A\times_\alpha G$. Let $\{P_\alpha(s);
s\in\widehat{G}\}$ be the spectral projections for $\lambda$:
$$
\lambda_t=\sum_{s\in\widehat{G}} \langle t,s\rangle P_\alpha(s)\;.
$$
Then it follows that $\iota_\ast([P_\alpha(0)])=[\1]$, where $\1$ is the
identity element of $\A$. Thus to each action $\alpha$ of $G$ on
$\A$ we associate the following $K$-theoretic data:
\begin{eqnarray*}
K_0(\A\times_\alpha G)
& \stackrel{\iota_\ast}{\longrightarrow}& K_0(\A) \\
\mbox{$[P_\alpha(0)]$} & \mapsto & [\1] \\
K_1(\A\times_\alpha G)
& \stackrel{\iota_\ast}{\longrightarrow} & K_1(\A)
\end{eqnarray*}
with the action $\widehat{\alpha}_\ast$ of $\widehat{G}$ on
$K_i(\A\times_\alpha G)$ satisfying
$\iota_\ast\widehat{\alpha}_\ast(s)=\iota_\ast$.

\begin{proposition}
Let $\A$ be a unital C*-algebra and let $\alpha,\beta$ be
actions of a compact abelian group $G$ on $\A$. If
$\alpha$ is conjugate to $\beta$, then there exist isomorphisms
$\varphi_1$ of $K_i(\A\times_\alpha G)$ onto
$K_i(\A\times_\beta G)$ and $\varphi_2$ of $K_i(\A)$ onto
$K_i(\A)$ such that
\begin{eqnarray*}
&&\iota_\ast\varphi_1=\varphi_2 \iota_\ast \\
&&\varphi_1\widehat{\alpha}_\ast(s)=
            \widehat{\beta}_\ast(s)\varphi_1 \\
&&\varphi_1([P_\alpha(0)])=[P_\beta(0)]
\end{eqnarray*}
where $\varphi_1,\varphi_2$ are order isomorphisms if applicable.
\end{proposition}

\begin{proof}
Note that $\varphi_2([\1])=[\1]$ follows automatically.

If $\alpha=\sigma\beta\sigma^{-1}$ for an automorphism $\sigma$ of
$\A$, then there is an isomorphism $\phi$ of
$\A\times_\alpha G$ onto $\A\times_\beta G$ such
that
\begin{eqnarray*}
\phi(a)=\sigma^{-1}(a),\qquad a\in\A \\
\phi(\lambda(t))=\lambda(t),\;\qquad t\in G
\end{eqnarray*}
By setting $\varphi_1=\phi_\ast$ and $\varphi_2=\sigma^{-1}_\ast$,
all the properties follow easily.
\end{proof}

\begin{proposition}
Let $\A$ be a unital C*-algebra and let $\alpha,\beta$ be
actions of a compact abelian group $G$ on $\A$. If
$\alpha$ is cocycle conjugate to $\beta$, then there exist
isomorphisms $\varphi_1$ of $K_i(\A\times_\alpha G)$ onto
$K_i(\A\times_\beta G)$ and $\varphi_2$ of
$K_i(\A)$ onto $K_i(\A)$ such that
\begin{eqnarray*}
&&\iota_\ast\varphi_1=\varphi_2 \iota_\ast \\
&&\varphi_1\widehat{\alpha}_\ast(s)=
    \widehat{\beta}_\ast(s)\varphi_1 \\
&&\varphi_2([1])=[1]
\end{eqnarray*}
where $\varphi_1,\varphi_2$ are order isomorphisms if applicable.
\end{proposition}

\begin{proof}
In view of Proposition~3.1 it suffices to show this when $\alpha$ is
a cocycle perturbation of $\beta$, i.e.,
$\alpha_t=\Ad u_t\circ\beta_t$ with $u_t$ a one-cocycle for $\beta$.

We define an isomorphism $\phi$ of $\A\times_\alpha G$
onto $\A\times_\beta G$ by
\begin{eqnarray*}
&&\phi(a)=a, \qquad a\in\A \\
&&\phi(\lambda(t))=u_t\lambda(t),\qquad t\in G
\end{eqnarray*}
Since $\widehat{\beta}(s)\phi=\phi\widehat{\alpha}(s)$,
$s\in\widehat{G}$, $\phi$ naturally extends to an isomorphism
$\widehat{\phi}$ of $\A\times_\alpha
G\times_{\widehat{\alpha}}\widehat{G}$ onto
$\A\times_\beta G\times_{\widehat{\beta}}\widehat{G}$. By
setting $\varphi_1=\phi_\ast$ and $\varphi_2=\widehat{\phi}_\ast$,
all the properties follow easily, perhaps except for
$\varphi_2([\1])=[\1]$. For this we shall show that $\varphi_2=\id$.

Supposing that $\A$ is represented on a Hilbert space
${\mathcal H}$, we represent $\A\times_\alpha
G\times_{\widehat{\alpha}} G$ on $L^2(G,{\mathcal H})$ by
\begin{eqnarray*}
&&(\pi_\alpha(a)\xi)(s)=\alpha_{-s}(a)\xi(s) \\
&&(\lambda(t)\xi)(s)=\xi(s-t) \\
&&(v(p)\xi)(s)=\langle p,s\rangle\xi(s)
\end{eqnarray*}
for $a\in\A$, $s,t\in G$, and $p\in\widehat{G}$. Then
$\A\times_\alpha G$ is generated by
$\pi_\alpha(\A)$ and $\lambda(f)$, $f\in L^1(G)$ and
$\A\times_\alpha G\times_{\widehat{\alpha}}\widehat{G}$ is
generated by $(\A\times_\alpha G)v(f)$, $f\in
L^1(\widehat{G})$, which naturally identifies with
$\A\otimes \mathfrak{K}(L^2(G))$. Then $\widehat{\phi}$ is given by
\begin{eqnarray*}
&& \pi_\alpha(a)\mapsto \pi_\beta(a) \\
&& \lambda(t)\mapsto u_t\lambda(t) \\
&& v(p)\mapsto v(p)
\end{eqnarray*}
and $\widehat{\phi}$ is implemented by the unitary $U$
defined by
$$
(U\xi)(s)=u_{-s}^\ast \xi(s)\;,
$$
which is a multiplier of
$\A\otimes \mathfrak{K}(L^2(G))$. Thus we
obtain the commutative
diagram:
$$
\xymatrix{
\A\times_\alpha G \ar[d]_\phi \quad
     \ar
@{^{(}->}^-{\iota}[r] & \quad \A\otimes \mathfrak{K}(L^2(G))
\ar[d]^{\widehat{\phi}=\Ad U} \\
\A\times_\alpha G\quad\ar
@{^{(}->}^-{\iota}[r] &
     \quad \A\otimes \mathfrak{K}(L^2)G))
}
$$
Thus the assertion follows.
\end{proof}

\section{Conjugacy and cocycle conjugacy classes of actions of $\T$}

As an application of Theorem~1.2 we consider the problem of
classifying a class of actions of $\T$ on a unital separable simple
purely infinite C*-algebra $\A$.

To meet the assumption in Theorem~1.2 we have to assume that
$\A\times_\alpha\T$ is AF, which implies in particular that $\A$ is
nuclear. By the Pimsner-Voiculescu exact sequence,
$\iota_\ast:K_0(\A\times_\alpha\T)\to K_0(\A)$ is
just the quotient map obtained by division by
$\Im(\id-\widehat{\alpha}_\ast)$, where we now denote by $\widehat{\alpha}$
the single automorphism
$\widehat{\alpha}(1)$ of $\A\times_\alpha\T$. Thus the
invariants for conjugacy classes described in 3.1 reduce to
$$
(K_0(\A\times_\alpha\T),[P_\alpha(0)],
\widehat{\alpha}_\ast)\;,
$$
and the following relations hold in this case:
\begin{eqnarray*}
&&(K_0(\A\times_\alpha\T)/\Im(\id-\widehat{\alpha}_\ast),\quad
     [P_\alpha(0)]+\Im(\id-\widehat{\alpha}_\ast))\cong
     (K_0(\A),[1]) \\ [.5ex]
&&\ker((\id-\widehat{\alpha}_\ast)|K_0(\A\times_\alpha\T)\cong
K_1(\A)\;.
\end{eqnarray*}

\begin{corollary}
Let $\A$ be a unital separable simple purely infinite
C*-algebra and let $\alpha,\beta$ be actions of $\T$ on $\A$ such
that $\A\times_\alpha\T$ and $\A\times_\beta\T$ are simple AF algebras
with one-dimensional lattices of traces. Then the following
conditions are equivalent:
\begin{enumerate}
\item
$\alpha$ and $\beta$ are conjugate

\item
$(K_0(\A\times_\alpha\T), [P_\alpha(0)],
\widehat{\alpha}_\ast)$ and $(K_0(\A\times_\beta\T),
[P_\beta(0)], \widehat{\beta}_\ast)$ are isomorphic, i.e., there
is an order-isomorphism $\varphi_1$ of
$K_0(\A\times_\alpha\T)$ onto
$K_0(\A\times_\beta\T)$ such that
$\varphi_1\widehat{\alpha}_\ast=\widehat{\beta}_\ast\varphi_1$ and
$\varphi_1([P_\alpha(0)])=[P_\beta(0)]$.
\end{enumerate}
\end{corollary}

\begin{proof}
We have shown (1)$\Rightarrow$(2) in 3.1.

Suppose (2). Note that $\A\times_\alpha\T$ and
$\A\times_\beta\T$ are stable AF algebras and
$\widehat{\alpha}_\ast,\widehat{\beta}_\ast,\varphi_1$ are all
isomorphisms such that
$$
\xymatrix{
K_0(\A\times_\alpha \T) \ar[d]_{\varphi_1} \quad
     \ar[r]^-{\widehat{\alpha}_\ast\;} &
\quad K_0(\A\times_\alpha\T)
       \ar[d]^{\varphi_1} \\
K_0(\A\times_\beta\T)\quad
     \ar[r]^-{\widehat{\beta}_\ast\;}
    & \quad K_0(\A\times_\beta\T)
}
$$
is commutative. Then by a standard intertwining argument we obtain
isomorphisms $\widehat{\alpha}',\widehat{\beta}',\phi$ such that
$$
\xymatrix{
\A\times_\alpha \T \ar[d]_{\phi} \quad
     \ar[r]^-{\widehat{\alpha}'\;}
    &  \quad \A\times_\alpha\T
       \ar[d]^{\phi} \\
\A\times_\beta\T\quad
     \ar[r]^-{\widehat{\beta}'\;}
    & \quad \A\times_\beta\T
}
$$
is commutative and $\widehat{\alpha}'=\widehat{\alpha}_\ast$,
$\widehat{\beta}'=\widehat{\beta}_\ast$, and $\phi_\ast=\varphi_1$.

Since $\A$ is purely infinite and
$\A\times_\alpha\T$ has only one trace $\tau$ up to
constant multiple, $\widehat{\alpha}$ must scale $\tau$, i.e.
$\tau\widehat{\alpha}=\lambda\tau$ with $\lambda\not=1$. The same is
true for $\widehat{\beta}$. Hence
by 1.2,
$\widehat{\alpha}'$,
$\widehat{\alpha}$, $\widehat{\beta}'$, $\widehat{\beta}$ have the
Rohlin property and hence $\widehat{\alpha}'$ and
$\widehat{\alpha}$ (resp. $\widehat{\beta}'$ and $\widehat{\beta}$)
are outer conjugate, i.e. $\widehat{\alpha}'=\Ad
u\sigma\widehat{\alpha}\sigma^{-1}$ and $\widehat{\beta}'=\Ad v
\nu\widehat{\beta}\nu^{-1}$ for unitaries $u, v$ in
$\A\times_\alpha\T+\1, \A\times_\beta\T+\1$ and automorphisms $\sigma,\nu$
of
$\A\times_\alpha\T , \, \A\times_\beta\T$, respectively. Moreover we can
assume that
$\sigma_\ast=\id$ and $\nu_\ast=\id$. Since
$\phi\widehat{\alpha}'=\widehat{\beta}'\phi$ we obtain
$$
(\Ad \phi(u))\phi\sigma\widehat{\alpha}\sigma^{-1}=(\Ad v)
\nu\beta\nu^{-1}\phi
$$
which implies, with $\phi_1=\nu^{-1}\phi\sigma$ and
$u_1=\nu^{-1}(\phi(u^\ast)v)$,
$$
\phi_1\widehat{\alpha}=\Ad u_1\widehat{\beta}\phi_1\;.
$$

Let $\{P_\alpha(n); n\in\Z\}$ (resp. $\{P_\beta(n); n\in\Z\}$) be the
spectral projections for $\lambda(t)$, $t\in\T$ in
$\A\times_\alpha\T$ (resp. $\A\times_\beta\T$), i.e.
$$
\lambda(t)=\sum_{n=-\infty}^\infty t^nP_\alpha(n)\;.
$$
Since $\widehat{\alpha}(P_\alpha(n))=P_\alpha(n-1)$ and
$\phi_{1\ast}([P_\alpha(0)])=[P_\beta(0)]$, we obtain that
$$
[\phi_1(P_\alpha(n))]=[P_\beta(n)]\;.
$$
Since $(\sum_{n=-N}^N P_\alpha(n))_N$ forms an approximate identity for
$\A\times_\alpha \T$, there is a unitary $U$ in $M(\A\times_\beta\T)$
such that
$\Ad U\circ\phi_1(P_\alpha(n))=P_\beta(n)$, $n\in\Z$. Thus replacing
$\phi_1$ by $(\Ad U)\phi_1$ we obtain an isomorphism $\phi$ of
$\A\times_\alpha\T$ onto $\A\times_\beta\T$ such
that
\begin{eqnarray*}
&&\phi(P_\alpha(n))=P_\beta(n),\qquad n\in\Z \\
&&\phi\widehat{\alpha}=(\Ad u)\widehat{\beta}\phi
\end{eqnarray*}
where $u$ is a unitary in $M(\A\times_\beta\T)$. Since
$$
P_\beta(n-1)=\phi\widehat{\alpha}(P_\alpha(n))=
\Ad u(P_\beta(n-1))\;,
$$
$u$ commutes with $P_\beta(n)$.
Let
$$
u_n=\left\{\!
\begin{array}{ll}u\widehat{\beta}(u)\ldots

\widehat{\beta}^{n-1}(u), & n=1,2,\ldots \\ [.5ex]
1 & n=0 \\
[.5ex]
\widehat{\beta}^{-|n|}(u_{|n|}^\ast)=

\widehat{\beta}^{-1}(u^\ast)\ldots
      \widehat{\beta}^{-|n|}(u^\ast), &

n=-1,-2,\ldots\end{array}\right.
$$
Then $\{u_n; n\in\Z\}$ commutes with
$P_\beta(n)$, and satisfies
that
$$
u_{m+n}=u_m\widehat{\beta}^m(u_n)\;,\qquad m,n\in\Z\;.
$$
We define
a unitary $v\in M(\A\times_\beta\T)$ by
$$
v=\sum u_n
P_\beta(-n)\;.
$$
Then it follows that
$$
v\widehat{\beta}(v^\ast)=\sum
u_n\widehat{\beta}(u_{n-1}^\ast)
P_\beta(-n)=u\;.
$$
Hence we obtain,
replacing $(\Ad v^\ast)\phi$ by
$\phi$
$$
\phi\widehat{\alpha}=\widehat{\beta}\phi
$$
and
$$
\phi(P_\alpha(n))=P
_\beta(n)\;.
$$
Then the extension of $\phi$ to the multiplier algebras
maps
$\A$ onto $\A$ and
$$
\phi(\alpha_t(a))=\beta_t(\phi(a))\;\qquad
a\in\A\;.
$$
(cf. \cite{Ped} 7.8.8). Thus $\alpha$ is conjugate to
$\beta$.
\end{proof}

\begin{corollary}
Let $\A$ be a unital separable
simple purely infinite
C*-algebra and let $\alpha,\beta$ be actions of $\T$
on
$\A$ such that $\A\times_\alpha\T$ and
$\A\times_\beta\T$ are simple AF
algebras with
one-dimensional lattice of traces. Then the following
conditions
are equivalent:
\begin{enumerate}
\item
$\alpha$ and $\beta$ are
cocycle conjugate;

\item
$(K_0(\A\times_\alpha\T)$,
$[P_\alpha(0)]+\Im(\id-\widehat{\alpha}_\ast)$,
$\widehat{\alpha}_\ast)$ and
$(K_0(\A\times_\beta\T)$,
$[P_\beta(0)]+\Im(\id-\widehat{\beta}_\ast)$,
$\widehat{\beta}_\ast)$ are isomorphic, i.e., there is an order
isomorphism $\varphi_1$ of $K_0(\A\times_\alpha\T)$ onto
$K_0(\A\times_\beta\T)$ such that
$\varphi_1\widehat{\alpha}_\ast=\widehat{\beta}_\ast\varphi_1$ and
$\varphi_1([P_\alpha(0)])-[P_\beta(0)]\in
\Im(\id-\widehat{\beta}_\ast)$.
\end{enumerate}
\end{corollary}

\begin{proof}
We have shown (1)$\Rightarrow$(2) in 3.2.

Assume that (2) holds. As in the proof of 4.1 we obtain an isomorphism
$\phi$ of $\A\times_\alpha\T$ onto $\A\times_\beta\T$ such that
$$
\phi\widehat{\alpha}=(\Ad u)\widehat{\beta}\phi
$$
for some unitary $u\in\A\times_\beta\T+\1$.

Since $(\A\times_\beta\T\times_{(\Ad u)\widehat{\beta}}\Z$, $\T , \,
((\Ad u) \widehat{\beta})^{\widehat{}})$ is isomorphic to
$(\A\times_\beta\T\times_{\widehat{\beta}}\Z$, $\T,
\widehat{\widehat{\beta}}\,)$,
which is again isomorphic to $(\A\otimes \mathfrak{K}(L^2(\T)),\T,
\beta\otimes\Ad\widetilde{\lambda})$, we obtain an isomorphism
$\widehat{\phi}$ of $\A\otimes \mathfrak{K}(L^2(\T))$ onto
$\A\otimes \mathfrak{K}(L^2(\T))$ such that
$\widehat{\phi}(\alpha_t\otimes\Ad\widetilde{\lambda}_t)=
(\beta_t\otimes\Ad\widetilde{\lambda}_t)\widehat{\phi}$, extending
$\phi:\A\times_\alpha\T\to\A\times_\beta\T$. The image of
$[P_\alpha(0)]\in K_0(\A\times_\alpha\T)$ in $K_0(\A\otimes
\mathfrak{K}(L^2(\T))=K_0(\A)$ is $[\1]$; and $\widehat{\phi}_\ast$
preserves this class, i.e., $[\widehat{\phi}(\1\otimes p)]=[\1\otimes p]$
where $p$ is a minimal projection in $\mathfrak{K}(L^2(\T))$.

Thus we find a unitary multiplier $U$ of $\A\otimes \mathfrak{K}(L^2(\T))$
such that $(\Ad U\widehat{\phi})$ is the identity on
$\1\otimes\mathfrak{K}(L^2(\T))$.
Then there is an automorphism $\sigma$ of $\A$ such that $(\Ad
U)\widehat{\phi}=\sigma\otimes\id$, from which it follows that
$$
\sigma\alpha_t\otimes\id=
\Ad U(\beta_t\otimes\Ad\widetilde{\lambda}_t)(U^\ast)
\beta_t\sigma\otimes\id\;.
$$
Using $u_t=U(\beta_t\otimes\Ad\widetilde{\lambda}_t)(U^\ast)\in\A$,
which is a one-cocycle for $\beta$, we obtain that
$$
\sigma\alpha_t=\Ad u_t\beta_t\sigma\;.
$$
This shows that $\alpha$ is cocycle conjugate to $\beta$.
\end{proof}

If $\A$ is not unital we still have the following:

\begin{corollary}
Let $\A$ be a separable simple purely infinite C*-algebra and let
$\alpha,\beta$ be actions of $\T$ on $\A$ such that
$\A\times_\alpha\T$ and $\A\times_\beta\T$ are simple AF algebras
with one-dimensional lattice of traces. Then the following
conditions are equivalent:
\begin{enumerate}
\item
$(\A\otimes \mathfrak{K}(L^2(\T)),\T,\alpha\otimes \id)$ is isomorphic to
$(\A\otimes \mathfrak{K}(L^2(\T)),\T,\beta\otimes \id)$

\item
$(\A\otimes \mathfrak{K}(L^2(\T))$, $\T,\alpha\otimes
\Ad\widetilde{\lambda})$  is isomorphic to $(\A\otimes
\mathfrak{K}(L^2(\T))$, $\T,\beta\otimes\Ad\widetilde{\lambda})$

\item
$(K_0(\A\times_\alpha\T),\widehat{\alpha}_\ast)$ is isomorphic to
$(K_0(\A\times_\beta\T),\widehat{\beta}_\ast)$.
\end{enumerate}
\end{corollary}

\begin{proof}
By general theory we obtain that (1)$\Rightarrow$(3) and
(2)$\Rightarrow$(3). From a part of the proof of 4.2 we obtain
(3)$\Rightarrow$(2).

Assume that (3) holds. Denoting by $\mathfrak{K}$ the compact operators on
$L^2(\T)$, we consider the systems $(\A\otimes
\mathfrak{K},\T,\alpha\otimes \id)$ and
$(\A\otimes \mathfrak{K},\T,\beta\otimes \id)$. Since $\A^\alpha$ is a
hereditary C*-subalgebra of $\A\times_\alpha\T$ and
$\A\times_\alpha\T$ is a stable simple AF algebra, we have that
$\A^\alpha\otimes \mathfrak{K}\cong \A\times_\alpha\T$. In the spectral
subspace
$$
\{x\in M(\A\otimes \mathfrak{K}): (\alpha_z\otimes \id)(x)=zx, z\in\T\}
$$
we find a unitary $U$ \cite{KT}. Note that $\A\otimes \mathfrak{K}$ is
generated by $\A^\alpha\otimes \mathfrak{K}$ and $U\A^\alpha\otimes
\mathfrak{K}$. Define an automorphism $\gamma$ of $\A^\alpha\otimes
\mathfrak{K}$ by
$\gamma=\Ad U$.

Then, as
$\widehat{\gamma}=\alpha\otimes \id, (\A^\alpha\otimes
\mathfrak{K}\times_\gamma\Z\times_{\widehat{\gamma}}\T,
\widehat{\widehat{\gamma}}\,)$ is isomorphic to
$((\A\times_\alpha\T)\otimes \mathfrak{K},
\widehat{\alpha}\otimes \id)$, which implies that
$(K_0(\A^\alpha\otimes \mathfrak{K}),\gamma_\ast)\cong
(K_0(\A\times_\alpha\T),\widehat{\alpha}_\ast)$.

In the same way we obtain a unitary $V$ in
$$
\{x\in M(\A\otimes \mathfrak{K}); \beta_z\otimes \id(x)=zx, z\in\T\}
$$
and that
$$
(K_0(\A^\beta\otimes \mathfrak{K}),(\Ad V)_\ast)\cong
(K_0(\A\times_\beta\T),\widehat{\beta}_\ast)\;.
$$
Hence, by using the Rohlin property for
$\Ad U|\A^\alpha\otimes \mathfrak{K}$
and $\Ad V|\A^\beta\otimes \mathfrak{K}$, we obtain an isomorphism
$\varphi$ of $\A^\alpha\otimes \mathfrak{K}$ onto $\A^\beta\otimes
\mathfrak{K}$ such that
$$
\varphi\circ\Ad U=\Ad vV\circ \varphi
$$
for some unitary $v\in\A^\beta\otimes \mathfrak{K}+\1$. Then $\varphi$
extends to an isomorphism $\widehat{\varphi}$ of $\A\otimes \mathfrak{K}$
onto $\A\otimes \mathfrak{K}$ such that
\begin{eqnarray*}
&&x\mapsto \varphi(x)\;, \\
&&Ux\mapsto vVx
\end{eqnarray*}
for $x\in\A^\alpha\otimes \mathfrak{K}$. It follows that
$\widehat{\varphi}(\alpha\otimes \id)=(\beta\otimes \id)\widehat{\varphi}$.
\end{proof}

A class of examples of actions of $\T$ on Cuntz algebras ${\mathcal O}_d$
(with
$d<\infty$) is obtained from quasi-free actions \cite{BJO}. Other
examples, which we consider in Example~4.4 and 4.5 below, are obtained
starting with trace-scaling automorphisms of AF algebras.

\begin{example}
For each $\lambda\in(0,1)\cup(1,\infty)$ we define an order on
$\Z[x,x^{-1}]$ by $p\geq q$ if
$p(\lambda,\lambda^{-1})>q(\lambda,\lambda^{-1})$ or $p\equiv q$.
Denote the ordered abelian group so obtained by
$\Z[x,x^{-1}]_\lambda$. Since $\Z[\lambda,\lambda^{-1}]$ is dense in $\R$,
$\Z[x,x^{-1}]_\lambda$ is a dimension group. Let
$\A_\lambda$ be the stable AF algebra whose dimension group is isomorphic to
$\Z[x,x^{-1}]_\lambda$ and $\alpha$ an automorphism of $\A_\lambda$
such that $\alpha_\ast$ corresponds to the multiplication by $x$.
Then $\A_\lambda\times_\alpha\Z$ is purely infinite \cite{R2} and
it is isomorphic to ${\mathcal O}_\infty\otimes \mathfrak{K}$ \cite{Phi}
since
$K_0(\A_\lambda\times_\alpha\Z)=\Z[x,x^{-1}]/(1-x)\Z[x,x^{-1}]
\cong\Z$.
By cutting down by a projection $p\in\A_\lambda$ with
$[p]=1\in\Z[x,x^{-1}]_\lambda$, we thus obtain an action
$\gamma_\lambda$ on
${\mathcal O}_\infty\cong p(\A_\Lambda\times_\alpha\Z)p$
from the dual action $\widehat{\alpha}$. The invariant for the
cocycle conjugacy class of
$({\mathcal O}_\infty,\T,\gamma_\lambda)$ is given by
$$
(\Z[x,x^{-1}]_\lambda,1+(1-x)\Z[x,x^{-1}]_\lambda,x)
$$
Note that each $\gamma_\lambda$ has a unique KMS state at the
inverse temperature $-\ln \lambda$.
\end{example}

\begin{example}
For $d=1,2,\ldots$ denote by $G$ the abelian group
$$
\Z[x,x^{-1}]+d\Z[x,x^{-1},(1-x)^{-1}]
$$
and for $\lambda\in(0,1)\cup(1,\infty)$, denote by $G_\lambda$ the ordered
abelian group $G$ with the order defined by: $p\geq0$ if
$p(\lambda,\lambda^{-1},(1-\lambda)^{-1})>0$ or $p\equiv 0$. Denote
by $\A_\lambda$ the stable AF algebra corresponding to $G_\lambda$
and by $\alpha$ an automorphism of $\A_\lambda$ corresponding to the
multiplication by $x$. Then $\alpha$ has the Rohlin property. Since
$$
G/(1-x)G\cong \Z/d\Z
$$
$\A_\lambda\times_\alpha\Z$ is isomorphic to ${\mathcal O}_{d+1}\otimes K$.
By cutting down $\A_\lambda\times_\alpha\Z$ by a projection
$p\in\A_\lambda$ with $[p]=1$ we obtain an action $\gamma_\lambda$
of $\T$ on ${\mathcal O}_{d+1}$ from the dual action $\widehat{\alpha}$.
Thus we obtain a one parameter family $({\mathcal
O}_{d+1},\T,\gamma_\lambda)$ of actions of $\T$ on ${\mathcal O}_{d+1}$
just for ${\mathcal O}_\infty$ in the
previous example.
\end{example}

\end{document}